\documentclass[a4paper,12pt]{amsart}
\usepackage{amsmath,amsfonts,amssymb}
\usepackage[all]{xy}
\usepackage{hyperref}
\usepackage[latin1]{inputenc}
\usepackage{psfrag}
\usepackage[dvips]{graphicx}

\usepackage{amsmath,amsfonts,latexsym}
\usepackage{amscd,amssymb,amsmath}
\usepackage{amsfonts}
\usepackage{amsbsy}
\usepackage{epsfig,afterpage}
\usepackage{psfrag}
\usepackage{graphicx}
\usepackage{colortbl}
\usepackage{cancel}
\usepackage[all]{xy}
\usepackage{hyperref}

\newcommand{\cdim}{{\rm coh.dim}}

\newtheorem{theorem}{\textbf{Theorem}}[section] 
\newtheorem{lemma}[theorem]{\textbf{Lemma}}     
\newtheorem{corollary}[theorem]{\textbf{Corollary}}
\newtheorem{proposition}[theorem]{\textbf{Proposition}}
\newtheorem{definition}[theorem]{Definition}

\newtheorem{remark}[theorem]{Remark}

\newcommand{\bea}{\begin{eqnarray*}}
\newcommand{\eea}{\end{eqnarray*}}
\newcommand{\zz}[1]{}

\newcommand{\bz}{\mathbb{Z}}

\newcommand{\br}{\mathbb R}
\newcommand{\bc}{\mathbb C}

\newcommand{\bt}{\mathbb T}
\newcommand{\bn}{\mathbb N}

\title[Bourgin-Yang versions of the Borsuk-Ulam theorem]
 {Bourgin-Yang versions of the Borsuk-Ulam theorem for  $p$-toral groups} %

\thanks{$^{1}$Supported  by the Polish Research Grants NCN 2011/03/B/ST1/04533 and NCN 2015/19/B/ST1/01458.}

\thanks{$^{2}$Supported by FAPESP of Brazil, Grants 2012/24454-8 and 2013/24845-0.}

\thanks{$^{3}$Supported by FAPESP of Brazil, Grants 2012/24454-8 and 2013/10353-8.}

\author[W. Marzantowicz, D. de Mattos and E. L.
dos Santos]{Wac{\l}aw Marzantowicz$^{1}$, Denise de Mattos$^{2}$,\\ Edivaldo L. dos Santos$^{3}$}

\address{Faculty of Mathematics and Computer Sci.,\, Adam Mickiewicz
University of Pozna{\'n},\, ul. Umultowska 87,\, 61-614 Pozna{\'n},
Poland.} \email{marzan@amu.edu.pl}
\address{Instituto de Ci\^encias Matem\'aticas e de Computa\c{c}\~ao, Universidade de S\~ao Paulo, Departamento de
Matem\'atica, Caixa Postal 668, S\~ao Carlos-SP, Brazil, 13560-970.}
\email{deniseml@icmc.usp.br} \address{Universidade Federal de S\~ao
Carlos, UFSCAR, Departamento de Matem\'atica, Caixa Postal 676,
S\~ao Carlos-SP, Brazil, 13565-905.} \email{edivaldo@dm.ufscar.br}


\begin{document}
\maketitle

\begin{abstract}
Let $V$ and $W$ be orthogonal representations of $G$ with $V^G=
W^G=\{0\}$. Let $S(V )$ be the sphere of $V$ and $f : S(V ) \to W$
be a $G$-equivariant mapping. We give an estimate for the dimension
of the set $Z_f=f^{-1}\{0\}$ in terms of $ \dim V$ and $\dim W$,  if
$G$ is the torus $\mathbb T^k$, or the $p$-torus $\mathbb Z_p^k$.
This extends the classical Bourgin-Yang theorem onto this class of
groups. \zz{Also we provide  a sufficient condition for the existence of
equivariant map between spheres of two orthogonal representations of
such a group.} Finally, we show that for any $p$-toral group $G$ and
a $G$-map $f:S(V) \to W$, with $\dim  V=\infty$ and $\dim W<\infty$,
we have $\dim Z_f= \infty$.
\end{abstract}


\section{Introduction}

Let $G$ be a compact Lie group and let  $V$, $W$ be two  orthogonal
representations of $G$ such that $V^G=W^G=\{0\}$ for the sets of
fixed points of $G$. Let $f: S(V) \to W$ be a $G$-equivariant
mapping. By $Z_f$, we denote the set $Z_f := \{v\in S(V)\, |\, f(v)=0\}.$

The problem of estimating the covering dimension of the set $Z_f$
was considered firstly  by C. T. Yang \cite{Ya1,Ya2} and
(independently) D. G. Bourgin \cite{Bourgin} for the case
$G=\mathbb{Z}_{2}$. Specifically, they proved that for a
$\bz_2$-equivariant mapping $f$ from the unit sphere $S(\br^m)$ in
$\br^m$ into $\br^n$, where the Euclidean spaces are considered as
representations of $\bz_2$ with the antipodal action,
\begin{eqnarray*}\dim Z_f \geq m-n-1,\end{eqnarray*}
where $\dim$ is the covering dimension. Consequently, it is a
strengthening of the classical Borsuk-Ulam theorem, which is its
direct consequence.

In \cite{Dold}, Dold extended the Bourgin-Yang problem to a
fibre-wise setting, giving an estimate for the set $Z_f =
f^{-1}\{0\}$, where $\pi:E\to B$ and $\pi^{\prime}:E^{\prime}\to B$
are vector bundles and $f : S(E)\subset E \to E^{\prime}$ is a
$\mathbb{Z}_{2}$-map, which preserve fibres $(\pi^{\prime}\circ
f=\pi)$. In \cite{Izydorek} and \cite{Nakaoka} this problem was
considered for the case of the cyclic group $G=\mathbb{Z}_{p}$ ($p$
prime), and in \cite{Denise1} for bundles $E\to B$ whose fibre has
the same cohomology (mod p) of a product of spheres. All these
results gave the  Bourgin-Yang  theorem  for $G=\mathbb{Z}_{p}$,
with $p$ prime,  if we take $B$ as a single point.

Recently, in \cite{MrMaSa} the authors  considered  the Bourgin-Yang
problem for the case that $G$ is a cyclic group of a prime power
order, $G = \bz_{p^k}$, $ k\geq 1$.   Based on the result of
\cite{Ba0} it was  proved \cite[Theorem 1.1]{MrMaSa} that if  $V, \,
W$ are two orthogonal representations of the cyclic group
$\bz_{p^k}$ and $f: S(V) \to W$ an equivariant map  then the
covering dimension $\dim(Z_{f}) = \dim(Z_{f}/G) \geq \phi(V, W )$,
where $\phi$ is a function depending on $\dim V$, $\dim W$ and the
orders of the orbits of actions on $S(V) $ and $S(W)$ (cf.
\cite[Theorems 3.6 and 3.9]{MrMaSa}). In particular, if $\dim W\,\,
< \,\, \dim V/p^{k-1}$, then $\phi(V,W)\geq 0$, which means that
there is no $G$-equivariant map from $S(V)$ into $S(W)$. Up to our
knowledge the above are all known classes of groups for which the
Bourgin-Yang theorem has been shown.

It is worth  pointing out that the classical Bourgin-Yang problem
studied here is similar  but different than  the Bourgin-Yang, or
correspondingly  Borsuk-Ulam problem for coincidence points along an
orbit of the action. The latter studied for $G=\bz_{p^k}$ by
Munkholm and for  $G=\bz_p^k$ by Volovikov in several papers (cf.
\cite{Munkholm1}, \cite{Munkholm2} and \cite{Vol1} with references there).   The
resulted outcomes are   called  the Borsuk-Ulam, or respectively
{\it Bourgin-Yang type} theorem.  It studies dimension of the set
$A(f)=\{x\in X:\, |\, f(x)=f(gx),\, {\text{for all}}\,\, g\in G\}$
for a map (not equivariant in general) $f: X\to Y$ between two
$G$-spaces $X$ and $Y$. There are relations between the estimations
of  dimensions derived for these two distinct problems but hopeless
to get each of them from the second directly  and we will not
discuss it.

 In this paper  we study the problem
of  estimation   the cohomological dimension of the set
$Z_f=\{x\in S(V): f(x)=0\}$ where $f:S(V) \to  W$ is an equivariant
map, and $V,W$ are orthogonal representations of the  group  $G$
such that $W^G=V^G=\{0\}$, where $G=\bz_{p}^k$ or $G=\bt^k$.

First in Section \ref{section2}, we prove in  Theorem \ref{B-Y for p-torus} that \[ \cdim  Z_f \geq \dim_{\br} V- \dim_{\br} W -1 .\,\]

 It gives an answer to the
classical Bourgin-Yang problem for this class of groups.

  As
an accompanying result we give a sufficient condition on $V$ and $W$
for the existence of equivariant map $f:S(V)\to S(W)$ which together
with earlier known classical necessary condition completes the
Borsuk-Ulam theorem for $p$-torus (Theorem \ref{necessary
condition for equivariant map Z_p}).

Finally, in Section \ref{section4} we discuss the Bourgin-Yang
problem for an action of a  $p$-toral group, i.e. a group $G$ of the
form $1\to \bt^k \to G \to \mathcal{P}\to 1 $, where $\mathcal{P}$
is a $p$-group.
 Theorem  \ref{characterization of
p-toral}  states that for any $G$-map $f:S(V) \to W$, of two
orthogonal representations of a $p$-toral group  with $\dim
V=\infty$ and $\dim W<\infty$, $V^G=W^G=\{0\}$ we have $\dim Z_f=
\infty$.   However here the argument is more complicated and uses
the mentioned result of \cite{BaClPu} based on the G.  Carlsson
theorem on the G. Segal conjecture and  theorem of Laitinen on the
completion of  Burnside ring of a $p$-group in the $p$-adic
topology.  In particular it  is purely infinite-dimensional, i.e.
does not give any estimate if $V$ is of finite dimension. Moreover,
combining it with the old result of T. Bartsch (\cite{Ba1}, and
\cite{Ba2}) one get a characterization of $p$-toral groups as the
unique class of groups with this property (Theorem
\ref{characterization of p-toral}).

We should  say that the Borsuk-Ulam theorem has many interesting
applications to the discrete mathematics (see  \cite{Matousek} for
details). Also one can deduce a generalization of the Tverberg
theorem (cf. \cite{Vol1}), or equipartition theorems as in
\cite{DolKar}. On the other hand there are many mini-max invariants
of a  $G$-spaces  which computations are based on the Borsuk-Ulam
theorem. They are used in several nonlinear variational  problems
with symmetry to estimate from below the number of solutions (cf.
\cite{Ba2} for a thorough exposition). Up to our knowledge the
Bourgin-Yang theorem does not have so spectacular  applications yet,
but there is an expectation for results giving estimates of the
dimension of the set of solutions of some problems where the
Borsuk-Ulam theorem gives only the existence of them.

Throughout the paper $\dim X$ stands for the covering dimension of a space $X$ and  $\cdim X$ stands for the cohomological dimension of a space $X$, i.e.,
$$\cdim X = \max\{n \,\,| \,\, \check{H}^{n}(X)\not= 0 \}$$
where $\check{H}^n(-)$ denotes the  \v{C}ech cohomology with coefficients $\mathbb{F} = \mathbb{Z}_p$ or $\mathbb{F} = \mathbb{Q}$, depending on whether $G=\bz_{p}^k$ or $G=\bt^k$.
 Since we are working with \v{C}ech cohomology theory, we have $\cdim X \leq \dim X$. Also, ${H}_*(-)$, ${H}^*(-)$ ($\tilde{H}_*(-)$, $\tilde{H}^*(-)$) denote the (reduced) singular (co)homology with
 coefficients $\mathbb{F} = \mathbb{Z}_p$ or $\mathbb{F} = \mathbb{Q}$, depending on whether $G=\bz_{p}^k$ or $G=\bt^k$.

 Next, let us consider the important Borsuk-Ulam type theorem proved by Assadi in \cite[page 23]{Assadi} (for $p$-torus) and Clapp and Puppe  in \cite[Theorem 6.4]{Clapp}.

\begin{theorem}\label{teorema2}
Let $G$ be a $p$-torus or a torus. Let $X$
and $Y$ be $G$-spaces with fixed-points-free actions; moreover, in
the case of a torus action assume additionally that $Y$ has finitely
many orbit types. Suppose that $\tilde{H}_j(X)=\tilde{H}^j(X) = 0$ for $j<n,$ $Y$
is compact or paracompact and finite-dimensional, and $H_j(Y)=H^j(Y)=0$ for
$j\geq n$. 
Then there exists no $G$-equivariant
map of $X$ into $Y$.
\end{theorem}

We recall that for $G=\bz_{p}^k$, with $p$ prime odd, and $G=\bt^k$ every
nontrivial irreducible orthogonal representation is even dimensional
and admits the complex structure, thus $V$ and $W$ admit it too.
Denote $d(V)=\dim_{\bc} V =\frac{1}{2} \dim_\br V$, and
corres\-pondingly $d(W)=\dim_{\bc} W =\frac{1}{2} \dim_\br W$. If
$G=\bz_{2}^k$ and $V$,$W$ are orthogonal representations of  $G$,
then denote $d(V)= \dim_\br V$, and respectively  $d(W)= \dim_\br
W$.

\section{\label{section2} Bourgin-Yang theorem for
$p$-torus  and torus}

The next result is the classical version of the Bourgin-Yang theorem for
$p$-torus and torus.

\begin{theorem}\label{B-Y for p-torus}
Let $V$, $W$ be two orthogonal representations of the group
$G=\bz_p^k$ or $G=\bt^k$ such that $V^G=W^G=\{0\}$. If  $f: S(V) \to W$ be a
$G$-equivariant map, then
\[ \cdim  Z_f \geq \dim_{\br} V- \dim_{\br} W -1 \,.\]

In
particular, if  $\dim_{\mathbb{R}} W \,\, < \,\, \dim_{\mathbb{R}}
V$, then there is no $G$-equivariant map from $S(V)$ into $S(W)$.
\end{theorem}

\begin{proof} Denote $m=\dim_{\br} V$ and $n=\dim_{\br} W$ and suppose
$$\cdim Z_{f}< m-n-1.$$

Then, $$\check{H}^i(Z_{f})=0, \,\,\text{for any}\,\, i > m-n-2.$$

By using Poincar\'e-Alexander-Lefschetz  duality and the long exact
sequence of the pair \linebreak $(SV, SV\setminus Z_{f})$, we
conclude
$$0=\check{H}^i(Z_{f})=H_{m-1-i}(SV, SV\setminus Z_{f})=\tilde{H}_{m-i-2}(SV\setminus Z_{f}), \,\,\text{for}\,\, j=m-i-2<n, \,\,\text{i.e.,}$$
$$\tilde{H}_{j}(SV\setminus Z_{f})=0, \,\,\text{for}\,\, j<n.$$

On the other hand, we have

$$H_{j}(W\setminus\{0\})=H_{j}(SW)=0, \,\,\text{for}\,\, j\geq n.$$

However, $$f:SV\setminus Z_{f}\to W\setminus \{0\}$$ is a
$G$-equivariant map, which contradicts Theorem \ref{teorema2}.

In particular, if $ \dim_{\mathbb{R}} V>
\dim_{\mathbb{R}} W$,  for a $G$-map $f:
S(V)\to S(W)\subset W$ it implies that $\cdim Z_f \geq
0$ and, consequently, $Z_f \neq \emptyset$, which gives a
contradiction.

\end{proof}

As a consequence we get the following corollary.

\begin{corollary}

Let $V$, $W$ be two orthogonal representations of the group
$G=\bz_p^k$ with $p>2$, or $G=\bt^k$, such that $V^G=W^G=\{0\}$. If
$f: S(V) \to W$ is a $G$-equivariant map and $\dim_{\br} V >
\dim_{\br} W $ then
\[ \cdim  Z_f \geq 1 \,.\]
\end{corollary}
\begin{proof} Indeed, since every nontrivial  orthogonal representation of $G=\bz_p^k$ with
$p>2$, or  $G=\bt^k$, has a complex structure, the integral number
$\dim_{\br} V- \dim_{\br} W -1 = 2(d(V)-d(W)) -1 $ is positive and
odd.
\end{proof}

Now, for the group $G=\bz_{p}^k$ and $H \subset G$ an isotropy
group, considering $f^{H} =  f_{|S(V^{H})} : S(V^{H}) \to W^{H}$ and
the set
$$Z_f^{H} := \{v\in S(V^{H})\, |\, f^{H}(v)=0 \}= Z_f \cap S(V^{H})$$
we have the following results.
\begin{corollary}\label{first main} Let $V, \, W$ be two orthogonal representations of
the group $G=\bz_{p}^k$, with $ V^G=W^G=\{0\}$, and let $f: S(V) \to
W$ be an equivariant map. Then, for every isotropy  group
$H \subset G$ of the action on $S(V)$, we have
$f(S(V^{H})) \subset W^{H}$, and for the cohomological
dimension
\[ \cdim(Z_{f}^{H}) \geq
\dim_{\mathbb{R}} V^{H}- \dim_{\mathbb{R}}
W^{H} -1. \]
\end{corollary}
\begin{proof} Let $H$ be an isotropy subgroup in $S(V)$ and $f^H:
S(V^H) \to W^H$ the restriction of $f$ to $S(V^H)$. Note that $V^H$
is a sub-representation of $G$, the action of $G$ on $S(V^H)$
factorizes through $K= G/H$. Here $H\simeq \bz_p^{k^\prime}$ and  consequently, $G/H \simeq
\bz_p^{\bar{k}}$ with $\bar{k}= k-k^\prime$. Moreover,
$f:S(V^H)\to W^H$ is a $K$-equivariant map and therefore, the result
 follows from Theorem \ref{B-Y for p-torus} applied to the 
group $\bz_p^{\bar{k}}$.
\end{proof}

We say that an isotropy subgroup $H\subset G$  is maximal if it
is maximal with respect to the inclusion. Note that if $G=\bz_{p}^k$
and the $G$-space is $S(V)$, $V^G=\{0\}$, then $H$ is maximal if and
only if it is $p$-subtorus of rank $k-1$. Indeed, these subgroups
are maximal and they appear as the isotropy subgroups, because for
an irreducible representation $V_\alpha \subset V$ given by
$\rho_\alpha: G \to \bz_p\subset \{z\in \bc: \vert z \vert =1\}$ and
each point $x\in S(V_\alpha)\subset S(V)$  we have $G_x = \ker
\rho_{\alpha}$.

\begin{theorem}\label{necessary condition for equivariant map Z_p}
Let $V, \, W$ be two orthogonal representations of
the group $G=\bz_{p}^k$, with $ V^G=W^G=\{0\}$. A necessary  and sufficient condition for the existence  of a
$\bz_p^k$-equivariant map $f: S(V) \to S(W)$ is
$$ \dim_{\mathbb{R}}  V^{H}  \leq \dim_{\mathbb{R}} W^{H}$$
for every maximal isotropy subgroup $H$ on $S(V)$.
\end{theorem}

\begin{proof} If  there is a $G$-map $f:S(V) \to S(W)\subset W$, then
$Z_f=\emptyset$, which gives $\dim Z_f^{H} = -1$, for every
$H$. By Corollary \ref{first main}, we have $\dim_{\mathbb{R}}
V^{H}  \leq \dim_{\mathbb{R}} W^{H}$, for every
maximal isotropy subgroup $H$,
which appears in the decomposition 
 of $V$.

The converse was already proved in \cite[Theorem 2.5]{Mr2}, in
another formulation. We present this proof. It is enough to
show that for every maximal subgroup $H\subset G$, under the
assumption $0 < \dim_{\mathbb{R}} V^H \leq \dim_{\mathbb{R}} W^{H}$,
there exists a $\bz_p$-map $f^H : S(V^H) \to S(W^H)$. Indeed, once
more, using the fact that the action of $G$ on $S(V^H)$ and $S(W^H)$
factorizes through $K= G/H \simeq \bz_p$, and any such a map $f^H$
is $G$-equivariant, it is sufficient to take the joint of maps
of the corresponding joints
$$f= {\underset{H} *}\,  f^H: S(V)= S({\underset{H}\oplus}\, V^H) = {\underset{H}
*} \,S(V^H) \; \to\; {\underset{H}
*} \,S(W^H)= S({\underset{H}\oplus}\, W^H) = S(W).$$

We have $V^H = \oplus_{j=1}^{p-1}\, l_j \, V_j$, with $l_j\geq 0$
 and $W^H = \oplus_{j=1}^{p-1}
\, \tilde{l}_j \, V_j$,  with  $\tilde{l}_j\geq 0$, where each $V_j$ is an irreducible representation of $G$. Since $ d(V^H)=
\sum_{j=1}^{p-1} \, l_j
 \leq \sum_{j=1}^{p-1} \, \tilde{l}_j = d(W^H)$,
it is enough  to show that for every $1 \leq j_1, j_2 \leq p-1$,
there exists a $\bz_p$-equivariant map from $S(V_{j_1}) \to
S(V_{j_2})$. Let $1 \leq j_1^{-1} \leq p-1 $ be the inverse of $j_1$
in $\bz_p^*$. It is easy to check that the map $  S^1 \ni z \mapsto
z^{j_1^{-1} j_2} \in S^1$ is  the required $\bz_p$-map.

\end{proof}

Now, we will consider  the problem of existence of $G$-equivariant maps
from the sphere $S(V)$ of an infinite dimensional representation $V$
into the sphere $S(W)$ of a finite dimensional representation $W$,
or the estimate of dimension of the set $Z_f$, for a $G$-equivariant
map  $f:S(V) \to W$.
\begin{theorem}\label{infinite for p-tori}
Let $V,\,W $ be an orthogonal representations of  a $p$-torus $G=
\bz_p^k$, $p$ prime, or the torus $G=\bt^k$, such that $V^G=\{0\}=
W^G$. If $\dim V = \infty$ and $\dim W < \infty$, then for every
$G$-equivariant map $f: S(V) \to W $ we have
\[ \dim \, Z_f \geq \cdim Z_f \,=\, \infty \,.\]
 In particular, there is no  $G$-equivariant map $S(V)\to S(W)$ under
this assumption.
\end{theorem}

\begin{proof} For a given $d\in \bn$, let us take a sub-representation
$V(d)\subset V$ such that $\dim_{\br} V(d)\geq d$. Restricting the
map  $f:S(V)\to W$ to $S(V(d))\subset S(V)$, we have an equivariant
map $ f_d : S(V(d))\to W$ with $Z_{f_d} \subset Z_f$. By Theorem
\ref{B-Y for p-torus}, $\dim Z_{f_d} \geq \cdim Z_{f_d} \geq d- \dim_\br W
-1$. Now using the monotonicity of  dimension we get

\[ \dim Z_f \geq \cdim Z_f \geq     \; {\underset{d\to \infty}\lim}\,  d-\dim W  -1   = \infty.\]
\end{proof}

\section{$p$-toral groups \label{section4}}

In this  section, we show that Theorem \ref{infinite for p-tori} 
can be extended on a larger class of groups called $p$-toral. The
main result will be formulated analogous to  \cite[Theorem
3.1]{Ba2}.

\begin{definition}\label{p-toral}\rm
A compact Lie group $G$ is called \textit{$p$-toral} if it is  of the form of
an extension
$$ 1\hookrightarrow  \bt^k \hookrightarrow G \to P \to 1 , $$
where $P$ is a finite $p$-group.
\end{definition}

In this section, we will use the $G$-index of $G$-spaces defined by
the Borel equivariant stable cohomotopy theory, i.e the theory
$ h_G^*(X,A) = \pi_s^*(X\times_G EG, A\times_G EG), $
where $\pi_s^*$ denotes the stable cohomotopy theory.
Following  \cite[5.4)]{Ba2}, as a family $\mathcal{B}$  of orbits
defining  value of this length index, we take
$\mathcal{B} = \{G/H: H\varsubsetneq G\}.$
Taking $I=h^*(pt)$, $h^*_G$ and $\mathcal{B}$ as above, the value of
the length index defined by the triple $\{\mathcal{B}, h_G^*, I\}$
at a pair of $G$-spaces $(X,X^\prime)$ will be denote by
$l(X,X^\prime)$.
\begin{theorem}[Characterization of $p$-toral groups]\label{characterization of p-toral}
{\phantom{ We have the following}} {\;\;\;}  \hskip 6cm \\
\vspace{-0.5cm}

\begin{itemize}
\item[\rm a)] {Let $G$ be a $p$-toral group $1 \hookrightarrow \bt^k  \to G \to
P\to 1$.  Then, for the sphere $S(V)$ of an infinite-dimensional
fixed point free  $G$-Hilbert space (orthogonal representation) $V$
and a finite dimensional orthogonal representation $W$ of $G$, such
that $W^G=\{0\}$, and a $G$-equivariant map $f: S(V) \to W$, we have
$$ \dim Z_f\; =\; l(Z_f)\;=\; \infty\,.$$}
\item[\rm b)] {If $G$ is
not $p$-toral, then there exist an infinite-dimensional fixed point
free $G$-Hilbert space $V$, a finite dimensional representation $W$
of $G$ with $W^G=\{0\}$ and an equivariant map $f: S(V) \to W$ such
that
$$ Z_f = \emptyset, \;\ {\text{e.g.}}\;\; \dim Z_f = -1 <
\infty\,.$$}
\end{itemize}
\end{theorem}
\proof The part b)  follows directly from \cite[Theorem 2)]{Ba1}. It
states that for any not $p$-toral group and every orthogonal Hilbert
representation $V$,  $V^G=\{0\}$, there exist an orthogonal
representation $W$, with $W^G=\{0\}$ and $\dim W < \infty$, and a
$G$-map $f: S(V) \to W\setminus\{0\}$. Therefore $Z_f = \emptyset$,
which proves part b).

 To show  a) we adapt the arguments of \cite{BaClPu} and
\cite{Ba1} exposed in an extended form in    \cite[Chapter 5]{Ba2}.
First, we have the following
\begin{proposition}[{\cite[5.11,
5.12]{Ba2}}]\label{infinitness of l(X)} For a $p$-group $P$ and a
contractible $P$-space $X$, we have $l(X)= \infty$.
\end{proposition} Note
that we do not require $X^P= \emptyset$. Recall that if $V$ is an
infinite-dimensional Hilbert space then  $X=S(V)$ is a metric
$G$-space which is contractible, because $S(V)$ is homeomorphic to
${\overset{\circ}D}(V)$. Consequently, it follows from Proposition
\ref{infinitness of l(X)} that  $l(S(V))= \infty.$

On the other hand, we have the following
\begin{proposition}[cf. {\cite[5.4]{Ba2}}]\label{finitness for sphere}
For every finite dimensional orthogonal representation $W$, with
$W^G=\{0\}$, and any $G$-length index as above, we have $l(S(W)) <
\infty$.
\end{proposition}

\begin{proof} Indeed, $S(W)$ is a compact $G$-space and each orbit of
$S(W)$ can be mapped into some element (orbit) of $\mathcal{B}$.
Now, the statement reduces to \cite[Corollary 4.9 b)]{Ba2}.
\end{proof}

Note that the statement in the last propositions  holds for every
equivariant cohomology theory $h^*$.

An essential step in our proof is the following result.

\begin{lemma}\label{finitnes of index for finite-dimensional}
Let $G$ be a finite group. If $X$ is a finite dimensional metric
$G$-space, with $X^G=\emptyset$, then $l(X)<\infty$.
\end{lemma}
\begin{proof}
Since $X$ is finite dimensional, so is $X/G$. Since $G$ is finite,
there is only a finite number of orbit types on $X$. Now, by the
Mostow theorem (cf. \cite[Theorem 10.1]{Br2}), there exist a finite
dimensional orthogonal representation $V$ of $G$ and a $G$-embedding
$\iota: X\to V$  of $X$ into $V$. Since $X^G=\emptyset$, we have
$\iota(X) \cap
 V^G= \emptyset$. Consequently, composing $\iota$ with $p_0^\perp: V \to V^G_\perp $, the orthogonal
projection of $V$  onto the orthogonal complement of $V^G$, we get
an equivariant map $\phi: X \to V^G_\perp\setminus \{0\}$. Now,
composing $\phi$ with the retraction $V^G_\perp \setminus \{0\} \to
S(V^G_\perp)$ we obtain a $G$-equivariant map $\psi : X\to
S(V^G_\perp)$.  Therefore, it follows from the monotonicity property
of the length index $l$ (cf. \cite[Theorem 4.7]{Ba2}) and
Proposition \ref{finitness for sphere} that $l(X)\leq
l(S(V^G_\perp))<\infty.$
\end{proof}

\begin{proof}[Proof of Theorem \ref{characterization of p-toral} a)] 
  Suppose first that $G$ is finite, i.e. $k=0$, and $G=P$ is a finite
$p$-group and let $f:S(V) \to W $ be a $P$-equivariant map. Since
$Z_f $ is closed $G$-invariant subspace  of $S(V)$, by continuity
property of the length index (cf. \cite[4.7 Continuity]{Ba2}), there
exists an open $P$-invariant neighborhood $Z_f \subset \mathcal{U}$
such that $l(Z_f)= l(\mathcal{U})$.

Denote $\mathcal{V}= S(V)\setminus Z_f$ which is a $P$-invariant
open subset of $S(V)$. Note that \linebreak $f: \mathcal{V} \to W\setminus
\{0\}$ is an equivariant map, and it follows from the monotonicity
property of the index that $l(\mathcal{V}) \leq l(W\setminus
\{0\})$. Also, by Proposition \ref{infinitness of l(X)}, we have
$l(S(V))=\infty$, and $l(S(W)) <\infty$, by Proposition
\ref{finitness for sphere} and the assumption $W^G=\{0\}$.

On the other hand $W\setminus \{0\}$ is $P$-equivariantly homotopy
equivalent to $S(W)$ and, consequently, has the same index. Using
the subadditivity property of the index, we get
$$ \infty \,= \,l(S(V))\, \leq\,  l(\mathcal{U}) + l(\mathcal{V})
\leq  l(Z_f) + l(S(W)).$$ Since $l(S(W)) < \infty$, we conclude that
$l(Z_f) = \infty\,.$  Note that $Z_f^G =\emptyset$, because
$S(V)^G=\emptyset$.  If $\dim Z_f < \infty$, then $l(Z_f) < \infty$
by Lemma \ref{finitnes of index for finite-dimensional}, which is a
contradiction. Thus, $\dim Z_f = \infty$.

Now, assume that $G$ is an extension $1 \to \bt^k \hookrightarrow G
\to P \to 1$, with $k\geq 1$. We distinguish two cases:

\vspace{0.2cm}

\centerline{either $\dim V^{\bt^k} =\infty \,,$ or  $ \dim V^{\bt^k}
< \infty \,.$}

\vspace{0.2cm}

First suppose  that $\dim V^{\bt^k} =\infty$. Note that $V^{\bt^k}$
has a natural action of $P= G/\bt^k$, with the fixed point set $V^G=
(V^{\bt^k})^P = \{0\}$. Moreover, the restriction
$f_{|S(V)^{\bt^k}}$ maps $S(V)^{\bt^k}$  $P$-equivariantly into
$W^{\bt^k} \subset W$. Applying  the previous  case for $G=P$, and
for the triple $(V^{\bt^k}, W^{\bt^k},\, f_{|S(V)^{\bt^k}})$  we
conclude that
$\dim Z_f \geq \dim Z_{f^{\bt^k}} = \infty.$

Now assume  that $\dim V^{\bt^k} <\infty\,.$ First observe that
$V^{\bt^k}$ is a $N(\bt^k)$ invariant subspace of $V$, where $N(H)$
is the normalizer of $H$ in $G$. But $\bt^k$ is a normal subgroup as
the component of the identity, thus $V^{\bt^k}$ is a
sub-representation of $V$. Let $V^\prime = V^{\bt^k}_\perp$ be the
orthogonal complement of $V^{\bt^k}$. By our assumption $\dim
V^\prime = \infty$.

 Following a standard argument also used in \cite[Proof of Theorem 3.1a]{Ba2}), we claim that
for any $p$-toral group $G$ and for every compact $G$-space $A$,
with $A^G=\emptyset$, there exists a finite $p$-group $P$ of $G$,
which acts on $A$ without fixed points: $A^P = \emptyset$. To see
this, observe that $G$ can be approximated by finite $p$-groups.
More precisely, for any natural number $s$ consider the set
$P^s:=\{g \in G:g^{p^s} =e\}$. \zz{Take $X= A \sqcup B$ the disjoint
union with a natural action of $G$.} If $p^s$ is a multiple of the
order of $G/G_0$, $G_0$ the component of the identity, then $P^s$ is
a subgroup of $G$, according to known results (see \cite[5.4]{Ba2}
for references). $P^s$ is obviously a finite $p$-group which is an
extension of $G_0 \cap  P^s$ by $G/G_0$, i.e. has the form $1 \to
P_s \hookrightarrow P^s \to P$ with $P_s = P^s \cap (G_0=\bt^k) $.
Moreover, it is clear that $A^{P^s} = \emptyset$, if $s$ is big
enough, because $A$ is a compact fixed point free $G$-space, i.e.
has only a finite number of isotropy orbit types.

Now, take $A= S(W)$ and $P^s$ as above. Then, $f: S(V^\prime ) \to
W$ is a $P^s$-equivariant map  and  $W^{P^s} = \{0\}$. Consider
$\tilde{V}=(V^\prime)^{P^s}$. If $\dim \tilde{V} =\infty$, then
$\dim Z_{f|S(\tilde{V})} = \dim S(\tilde{V}) =\infty$, because
$f(S(V)^{P^s})\subset W^{P^s}= \{0\}$. Consequently, assume that
$\dim \tilde{V} <\infty$.  Once more, $\tilde{V}$ is a
$P^s$-subrepresentation and we can take the orthogonal complement
$V^{\prime\prime}$ of it in $V^\prime$. By our assumption and the
choice of $V^{\prime\prime}$, we have $(S(V^{\prime\prime}))^{P^s}=
\emptyset$ and $\dim V^{\prime\prime}= \infty$.

In this way, we reduced the assumption to the already studied case
of a finite $p$-group. Applying it to $f:S(V^{\prime\prime}) \to W$,
we have $\dim Z_f \geq \dim Z_{f_{|S(V^{\prime\prime})}} = \infty,$
which completes the proof of Theorem \ref{characterization of
p-toral}. \end{proof}
\begin{remark}\rm One can easily note that our proof of Theorem \ref{characterization of
p-toral} a) is more complicated than the corresponding Borsuk-Ulam
theorem presented in \cite{BaClPu} and \cite{Ba2}. It is caused by
the fact that we need the assumption $S(V)^{P^s}= \emptyset$, which
is not necessary in the study of the Borsuk-Ulam problem. Indeed, if
$S(V)^{P^s} \neq \emptyset$, then there is no $P^s$-map from $S(V)$
into $S(W)$, because $S(W)^{P^s}=\emptyset$. The mentioned
assumption is necessary to know that $Z_f^{P^s} = \emptyset$, which
is necessary to apply Lemma \ref{finitnes of index for
finite-dimensional}.
\end{remark}

\begin{remark}\rm
Note that the torus $G= \bt^k$ and the $p$-torus $\bz_p^k$ are
toral, with $p=1$ or $k=0$, respectively. Then, the above argument
can be applied to any prime $p$. In this way, we obtain another
proof of Theorem \ref{infinite for p-tori}. On the other hand, our
proof of Theorem \ref{infinite for p-tori} is completely elementary.
Contrary to it, the proof of Theorem \ref{characterization of
p-toral} uses Proposition \ref{infinitness of l(X)}. The latter is
the result of \cite{BaClPu} and is based on a deep topological
result namely the Segal conjecture, proved by G. Carlsson.
\end{remark}

\begin{corollary}
The statement of Theorem \ref{characterization of p-toral} holds, if
we replace $\dim Z_f$ by $\dim Z_f/G$ in its statement.
\end{corollary}
\begin{proof} The statement of corollary follows directly from  the main
result of \cite{Deo}, since for a compact Lie group $G$ and a metric
$G$-space $X$, we have $ \dim \, X - \dim \, X/G \leq \dim\, G$.
\end{proof}

\end{document}